\documentclass{article}
\usepackage{graphicx} 
\usepackage{amsmath, amsfonts}

\title{On a Class of Integrals related to Wallis Integrals}
\author{Andie Burazin, Veselin Jungi\'c, Miroslav Lovri\'c}

\begin{document}

\maketitle

\begin{abstract}
   In this article we provide a motivation for studying integrals $\displaystyle \int_0^{\frac{\pi}{2}}x\cos^mx~dx$, $m\in \mathbb{N}$. We demonstrate several properties of this class of integrals and show that it is closely related to the class of Wallis integrals, $\displaystyle \int_0^{\frac{\pi}{2}}\cos^mx~dx$, $m\in \mathbb{N}$.
\end{abstract}
\section{Introduction}

In the theory of sound waves, a sawtooth wave is described as a wave that ramps upwards
and then sharply drops; see Figure 1 (left).
The problem we study here stems from
an investigation of
a modification of such a wave,
in which the ramp (rise) is circular, and the amplitudes are limited
by a specific line; see Figure 1 (right).

\begin{figure}[h!]
\centering
\includegraphics[width=90mm]{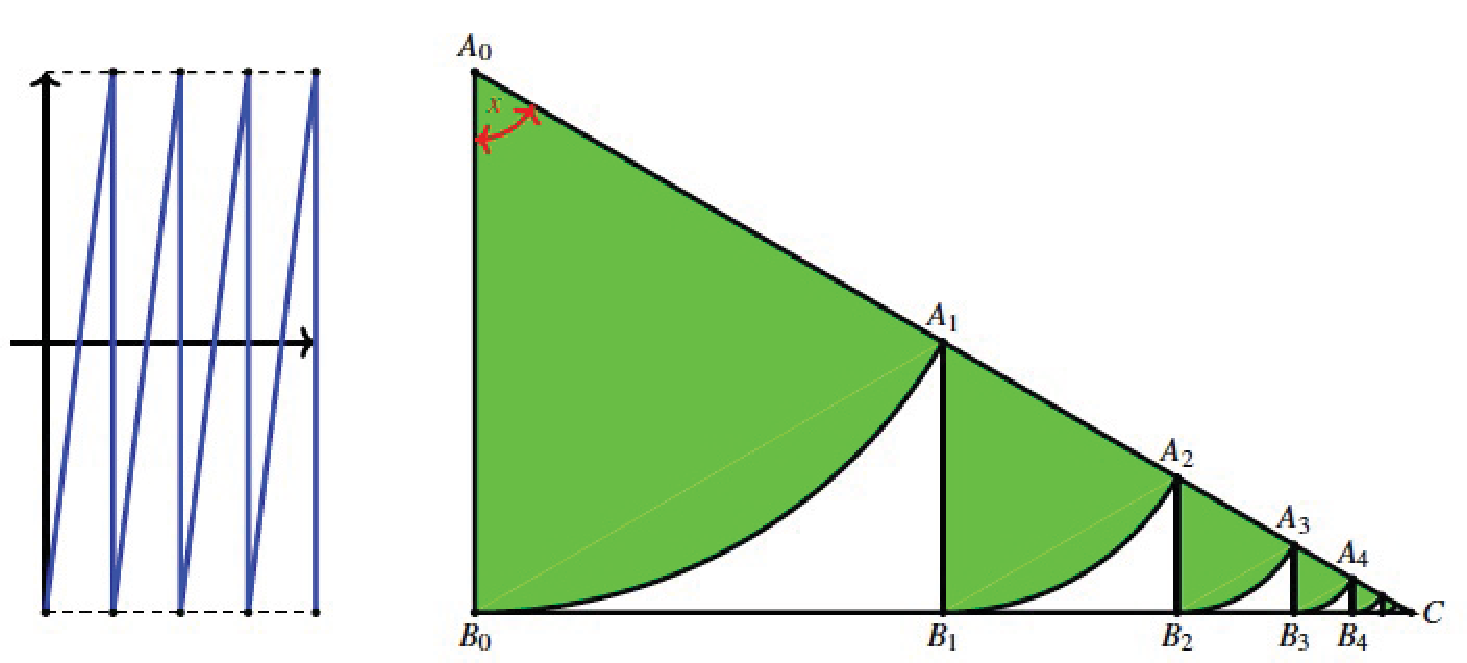}
\caption{
A sawtooth wave (left) and a modified sawtooth wave (right).}
\end{figure}

To construct this wave, we start
with a right triangle $\Delta A_0B_0C$
whose hypotenuse $\overline{A_0C}$  is of length 1, and the angle
$\angle CA_0B_0=x, $ where $0<x<\pi/2.$
For 
$i=1,2,3, \ldots$
we construct,
recursively,    
a sequence of points 
$A_i\in\overline{A_0C}$ and $B_i\in\overline{B_0C}$ so that
$|\overline{A_{i-1}A_i}|  =|\overline{A_{i-1}B_{i-1}}|$
and
$\overline{A_{i}B_{i}}   \perp \overline{B_0C}.$  

Among other quantities we studied, we looked into the areas of the circular sectors (Figure 1 (right), shaded).
It is a straightforward exercise to
show that the area of the $i$-th circular sector $A_iB_iA_{i+1}$ is  
$$  
A_i(x) =  {1 \over 2}\, x\cos^2 x (1-\cos x)^{2i}. 
$$  
To get a better feel about the function $x \mapsto A_i(x)$, for $x \in (0,\pi/2)$, we look at
its 
average value 
$$
\tilde{A_i} =   {1 \over \pi}\,  \int_0^{\pi/2} x\cos^2 x (1-\cos x)^{2i} \, dx .
$$
After expanding the integrand using Newton's binomial formula  
\begin{align}
\tilde{A_i} 
&
 = \frac{1}{\pi}\,  \int_0^{\pi/2}  x \cos^2 x    \sum_{m=0}^{2i} {2i \choose m} (-\cos x)^{m}\, dx  \nonumber \\ 
&   = \frac{1}{\pi}\,   \sum_{m=0}^{2i}(-1)^m {2i \choose m} \int_0^{\pi/2}     x  \cos ^{m+2} x \, dx  \label{jedan} 
\end{align}
we are left 
with a linear combination of integrals
of the form
$$ I_k = \int_0^{\pi/2} x \cos^kx \, dx$$
where $k\geq 2.$

Our original idea (discussed here as well)
was to prove that \eqref{jedan} can be written in the form
$$\tilde{A_i}  = \alpha_i +\beta_i \pi + { \gamma_i \over \pi}$$
where
$\alpha_i,$ $\beta_i, $ and $\gamma_i$ are rational numbers.

\medskip\medskip

\noindent
{\bf Problem.} 
What are the values of the integrals $I_k$? What are their properties?

\medskip
It turns out that we can say quite a lot, especially when $k$ is an even number.
(Although in this paper
we focus on even powers
of the cosine, we will make a few comments about odd powers
as well.)

\bigskip

Note that the {\it Wallis integrals}\footnote{There seems to be no agreement on whether to use the term \emph{Wallis integrals} or \emph{Wallis'  integrals}. We decided to use the former.}

$$ W_k = \int_0^{\pi/2}  \cos^kx \, dx$$
defined for $k \geq 0$ 
can be easily calculated 
from the recurrence formulas 
(see [3], section 7.1):

\begin{equation}
 W_k    
= \left\{ 
\begin{array}{lr}
 \displaystyle{{(k-1)!! \over k!!}{\pi \over 2}} & \hbox{ if $k$ is even, $k \geq 2$} \\
    \displaystyle{{(k-1)!! \over k!!}}  & \hbox{ if $k$ is odd, $k \geq 1$} \end{array}\right.  \label{dva}  
\end{equation}

Clearly, $W_0=\pi/2.$
See  [1] and [2] for related facts, 
including expressing $W_k$
in terms of the gamma function.

\section{Calculating $I_k$}

First, we derive a useful iterative formula for the definite integrals $I_k.$

\medskip

We start with the expression 
$$ \left(\sin x \cos^{k-1}x\right)' = k \cos^k x - (k-1)\cos^{k-2}x,$$
multiply it by $x,$ and rearrange terms:
$$ x\cos^k x  = {1 \over k} x \left(\sin x \cos^{k-1}x\right)'  + {k-1 \over k} x\cos^{k-2}x$$
Integrating both sides, we obtain
$$I_k = \int_0^{\pi/2} x \cos^kx \, dx
    = {1 \over k} \int_0^{\pi/2}  x \left(\sin x \cos^{k-1}x\right)' \, dx +  {k-1 \over k} I_{k-2}$$
To evaluate the remaining integral,
we use integration by parts:

\begin{align}
\int_0^{\pi/2}  x    \left(\sin x \cos^{k-1}x\right)' \, dx 
 & = \Bigl.  \sin x \cos^{k-1}x   \Bigr|_0^{\pi/2}   
           -  \int_0^{\pi/2} \sin x \cos^{k-1}x  \,dx   \nonumber \\ 
& = 0 + \left.{1 \over k} \cos^{k}x\right|_0^{\pi/2} 
  = - {1 \over k}.  \nonumber 
\end{align}

\noindent
Therefore, when $k \geq 2,$
\begin{equation}
  I_k =  {k-1 \over k} I_{k-2} - {1 \over k^2}  \label{tri}
\end{equation}

\noindent
Note that the well-known recurrence relation for
Wallis integrals, $$W_k={k-1\over k} W_{k-2},$$ for $k \geq 2,$
may be obtained in a similar way.

\bigskip
Recurrence \eqref{tri}, together with
$$ I_0 = \int_0^{\pi/2} x \, dx = {1 \over 8} \pi^2$$ 
and
$$ I_1 = \int_0^{\pi/2} x \cos x \, dx = 
   \Bigl. x \sin x + \cos x \Bigr|_0^{\pi/2}  
= -1 + {1 \over 2} \pi$$
gives
\begin{align}
    I_2 & = \int_0^{\pi/2} x \cos^2  x \, dx = {1\over2}I_0 - {1 \over 4} = - {1 \over 4} + {1 \over 16}\pi^2 \cr
I_3 & = \int_0^{\pi/2} x \cos^3  x \, dx = {2\over3}I_1 - {1 \over 9} = - {7 \over 9} + {1 \over 3}\pi  \cr
I_4 & = \int_0^{\pi/2} x \cos^4  x \, dx = {3\over4}I_2 - {1 \over 16} = - {1 \over 4} + {3 \over 64}\pi^2  \cr
I_5 & = \int_0^{\pi/2} x \cos^5  x \, dx = {4\over5}I_3 - {1 \over 25} = - {149 \over 225} + {4 \over 15}\pi  \nonumber
\end{align}
and so on.
In general,
from \eqref{tri} and the values of $I_0$ and $I_1$ 
we conclude that, for an even number $k\geq 0,$  
\begin{equation}
 I_k = \int_0^{\pi/2} x \cos^kx \, dx  = a_k\pi^2+ b_k  \label{cetri}
\end{equation}
where $a_k$ and $b_k$ are rational numbers, and 
\begin{equation}
I_k  = \int_0^{\pi/2} x \cos^kx \, dx  = c_k\pi + d_k   \label{oet}
\end{equation}
for and odd number $k\geq 1,$  where $c_k$ and $d_k$ are rational numbers.

\bigskip

The recurrence relation \eqref{tri} allows us to calculate $a_k,$ $b_k,$ $c_k$ and $d_k.$    
Using induction, we obtain, for $n \geq 1$ and $k \geq 2n,$

\begin{equation}
I_k  = I_{k-2n} \cdot \prod_{j=0}^{n-1} {k-2j-1 \over k-2j}
  - \left( {1 \over k^2} +\sum_{i=1}^{n-1} {1 \over (k-2i)^2} \prod_{j=0}^{i-1} {k-2j-1 \over k-2j} \right)    
  \label{sest} 
\end{equation}  
adopting the convention (and using it, from now on) that $\sum_{i=1}^{0}$  is zero.
 
\bigskip

When $k=2n,$ formula \eqref{sest} gives
\begin{equation}
I_{2n} = I_{0} \cdot \prod_{j=0}^{n-1} {2n-2j-1 \over 2n-2j}
  - \left( {1 \over 4n^2} +\sum_{i=1}^{n-1} {1 \over (2n-2i)^2} \prod_{j=0}^{i-1} {2n-2j-1 \over 2n-2j} \right)  
  \label{sedam}
\end{equation}
with $n \geq 1.$ 
For instance, using \eqref{sedam} with $n=1$ we obtain
$$
I_{2} =  I_0 \cdot {1 \over 2}  -  {1 \over 4} 
  = {1 \over 16}\pi^2 - {1 \over 4}
$$
which agrees, of course, with the value we computed earlier.

\bigskip

\noindent
Simplifying the first product in \eqref{sedam}, we obtain
$$
  \prod_{j=0}^{n-1} {2n-2j-1 \over 2n-2j}  =  
      {2n-1 \over 2n}\cdot {2n-3 \over 2n-2}\cdot{2n-5 \over 2n-4}\cdot\cdots\cdot{1 \over 2} 
         = {(2n-1)!! \over (2n)!!} 
$$
and therefore (recall that $ I_0 =  \pi^2/8$)
\begin{equation}
a_{2n}   
  = {1 \over 8}  {(2n-1)!! \over (2n)!!} 
 = {1 \over 4\pi} W_{2n}
   \label{osam}
\end{equation}
for $n \geq 1.$

\bigskip

The second product in formula \eqref{sedam} for $I_{2n}$ can be simplified in a similar way: 
expand the product, and multiply and divide by $(2(n-i)-1)!!$ and $(2(n-i))!!$ 
to obtain 
$$
\prod_{j=0}^{i-1} {2n-2j-1 \over 2n-2j}  = 
  {(2n-1)!! \,(2(n-i))!! \over (2n)!!(2(n-i)-1)!!} 
$$
and thus, when $n \geq 1,$

\begin{align}
b_{2n}
&
 =   - \left( {1 \over 4n^2} 
  +    {(2n-1)!! \over 4(2n)!!} \sum_{i=1}^{n-1} { 1\over (n-i)^2} { \,(2(n-i))!! \over (2(n-i)-1)!!} \right) \nonumber \\ 
&=  - \left( {1 \over 4n^2} 
  +    {(2n-1)!! \over 4(2n)!!} \sum_{j=1}^{n-1} { 1\over j^2} { \,(2j)!! \over (2j-1)!!} \right) \nonumber \\
  & =  -  {(2n-1)!! \over 4(2n)!!} \sum_{j=1}^{n} { 1\over j^2} { \,(2j)!! \over (2j-1)!!} \nonumber
\end{align}

\noindent
Because $I_0=\pi^2/8,$ it follows that $b_0=0.$

\bigskip

Thus, we have shown that,
for $k=2n\geq 2,$  
\begin{align}
I_{2n} & = \int_0^{\pi/2} x \cos^{2n}x \, dx  = a_{2n}\pi^2+ b_{2n} \nonumber \\
    & = {1 \over 8}  {(2n-1)!! \over (2n)!!} \pi^2
  -  {(2n-1)!! \over 4\cdot (2n)!!} \sum_{j=1}^{n} { 1\over j^2} { \,(2j)!! \over (2j-1)!!}  \nonumber \\
     & =   {\pi \over 4} W_{2n}
  -  {1 \over 2\pi} W_{2n} \sum_{j=1}^{n} { 1\over j^2} { \,(2j)!! \over (2j-1)!!} \nonumber \\
     & =   {1\over 4} W_{2n}
     \left( \pi  
  -  {2 \over \pi}   \sum_{j=1}^{n} { 1\over j^2} { \,(2j)!! \over (2j-1)!!}
   \right)
\label{devet}
\end{align}

\noindent
The fraction involving double factorials can be simplified:
$$
   {(2j)!! \over (2j-1)!!} 
     = { 2^j \cdot j(j-1)(j-2)\cdot\cdots\cdot 1 \over (2j-1)(2j-3)(2j-5)  \cdot\cdots\cdot 1} 
      \cdot {(2j)!! \over (2j)!!}  
     = {2^{2j} \, j! \, j! \over (2j)!} = {2^{2j} \over {2j \choose j}}       
$$

The formula for $I_k$ when $k$ is an odd number is derived similarly.
As we focus on even values of $k,$
we will not demonstrate it here.

\section{Investigating $I_k$}

We now
relate definite integrals
$$ I_k = \int_0^{\pi/2} x \cos^kx \, dx
\quad \hbox{and} \quad
J_k = \int_0^{\pi/2} x \cos^k 2x \, dx.$$
(Note that $I_0=J_0=\pi^2/8.$)

\medskip

Using the substitution $u=2x,$ 
$J_k$ can be written as
$$
J_k  = {1 \over 4} \int_0^{\pi} u \cos^ku \, du 
    = {1 \over 4}  \left(  I_k+ \int_{\pi/2}^{\pi} u \cos^ku \, du \right)
 $$
Now substitute $u=\pi-x$ to obtain
\begin{align}
J_k & =  {1 \over 4}  \left( I_k+   \int_0^{\pi/2} (\pi-x) (-1)^k \cos^kx \, dx   \right)  \nonumber \\
    & =  {1 \over 4}  
    \left( I_k - (-1)^k \int_0^{\pi/2} x \cos^kx \, dx  +(-1)^k\pi \int_0^{\pi/2} \cos^kx \, dx   \right)  \nonumber 
\end{align}
Thus,
\begin{align}
J_k & = {1 \over 4} \left( I_k +(-1)^{k+1}I_k + (-1)^{k} \pi \int_0^{\pi/2}  \cos^kx \, dx \right)  \nonumber \\
& = {1 \over 4} \Bigl( I_k +(-1)^{k+1}I_k + (-1)^{k} \pi W_k \Bigr). 
\label{deset}
\end{align}

\noindent
In particular, when $k\geq2$ is even, using \eqref{deset}, we obtain
\begin{equation}
J_{k}  = {\pi \over 4}W_k
    = {(k-1)!! \over k!!}\,{\pi^2 \over 8} 
\label{jedanaest}
\end{equation}
and, of course, 
$J_0  = {\pi^2 / 8}.$

\medskip

Thus, when $k\geq2$ is even, we write \eqref{devet} as 
\begin{equation}
I_{k}  = \int_0^{\pi/2} x \cos^kx \, dx  
     = J_{k}
  -  {(k-1)!! \over 4\cdot k!!} \sum_{j=1}^{k/2} { 1\over j^2} { \,(2j)!! \over (2j-1)!!}
\label{dvanaest}
\end{equation}
For instance 
(we know that $I_0 = J_0 ={\pi^2 / 8} $) 
using \eqref{dvanaest} and \eqref{jedanaest},
\begin{align}
I_2 & = J_2  - {1 \over 4} =  {1 \over 16} \pi^2  - {1 \over 4}  \nonumber \\
I_4 & =    J_4 - {1 \over 4}  =   {3 \over 64} \pi^2 - {1 \over 4} \nonumber \\
I_6 & = J_6 - {17 \over 72} =  {5 \over 128}  \pi^2 - {17 \over 72} \nonumber \\
I_8 & = J_8- {2 \over 9} =  {35 \over 1024} \pi^2 - {2 \over 9} \nonumber 
\end{align}
and so on.

\medskip
By using formula \eqref{jedanaest}, we  can rewrite \eqref{dvanaest} as:
$$
I_{k} 
     = J_{k}
  -  {1\over 4}\cdot {8J_k\over \pi^2} \sum_{j=1}^{k/2} { 1\over j^2} { \,(2j)!! \over (2j-1)!!}
     = J_{k}
 \left( 1 - {2\over \pi^2} \sum_{j=1}^{k/2}  {2^{2j} \over j^2 {2j \choose j}}  \right)   
$$

\noindent
Note: from \eqref{jedanaest} it follows that, due to symmetry, 
$$
J_{k}  = {\pi \over 4}\int_0^{\pi/2}  \cos^kx \, dx 
  = {\pi \over 4}\cdot{1\over 4}\int_0^{2\pi}  \cos^kx \, dx   
  = {\pi^2\over 8}\cdot {1\over 2 \pi}\int_0^{2\pi}  \cos^kx \, dx   
$$
i.e., $J_k$ is equal to the average value of $f(x)=\cos^kx$ on $[0,2\pi]$
multiplied by $\pi^2/8.$

\bigskip

We mention that when $k$ is odd, \eqref{deset} gives a relation between $I_k$ and $J_k$ in that case:
\begin{equation}
J_k = {1 \over 4} \left( 2I_k - \pi W_k \right)  =
  {1 \over 2 }I_k - {(k-1)!! \over k!!}\,{\pi  \over 4} 
\label{trinaest}
\end{equation}

\section{Several interesting results}

In this section we assume that $k \geq 2$ is even.

\medskip\medskip

\noindent
({\it i}) Write  
$$
  \cos^k x  = 
  \left( \cos^2 x \right)^{k/2}
  = \left( {1+\cos 2x \over 2} \right)^{k/2} = {1 \over 2^{k/2} } \left( 1+\cos 2x \right)^{k/2}
$$
and expand using the binomial formula, to
obtain the following 
relations:
\begin{align}
I_2 & = \int_0^{\pi/2} x \cos^2  x \, dx = {1 \over 2} \int_0^{\pi/2} x (1+\cos 2x) \, dx  
  = {1\over 2} (J_0+J_1) \nonumber \\
I_4 & = \int_0^{\pi/2} x \cos^4  x \, dx = {1 \over 2^2} \int_0^{\pi/2} x (1+\cos 2x)^2 \, dx  
  = {1\over 4} (J_0+2J_1+J_2)\nonumber \\
I_6 & = \int_0^{\pi/2} x \cos^6  x \, dx =   {1\over 8} (J_0+3J_1+3J_2+J_3) \nonumber \\  
I_8 & = \int_0^{\pi/2} x \cos^8  x \, dx = {1\over 16} (J_0+4J_1+6J_2+4J_3+J_4) \nonumber  
\end{align}
Note the naturally emerging Pascal's triangle pattern. In general, 
$$
I_{2n}  = \int_0^{\pi/2} x \cos^{2n}  x \, dx = {1\over 2^n} \sum_{i=0}^{n} {n \choose i} J_i  
$$
with $n \geq 1.$

\bigskip

\noindent
({\it ii}) 
Looking at the values of $J_k$ in \eqref{jedanaest}, we
remember that we
have seen these numbers before. It turns out
that
$J_k,$ for $k$ even, are the coefficients in the Taylor expansion: 
\begin{equation}
{1 \over \sqrt{1-x}}  = {8 \over \pi^2} \left(J_0+J_2x+J_4x^2+J_6x^3+\cdots\right)  \label{cetrnaest}
\end{equation}
Although this
is straightforward to check algebraically 
(see [2], where they establish the relation 
$$ {1 \over \sqrt{1-x}} = \sum_{j=0}^{\infty} {(2n-1)!! \over (2n)!!} x^n$$
in equation (23) on page 186) 
we do not know why, intuitively, something like this would work.

\smallskip
Connecting with the binomial formula, we get that
for $k=2n$ even,
$$ J_{2n} = (-1)^{n} {-1/2 \choose n} {\pi^2 \over 8}$$
As well, 
recalling the Maclaurin series formula,
\eqref{cetrnaest} gives 
$$ J_{2n} = \left. {d \over dx^n} \left( {1 \over \sqrt{1-x}} \right) \right|_{x=0}$$

\bigskip

\noindent
({\it iii}) 
When $k=2n\geq0$  is even, using \eqref{dva}, we note that
$$
W_{2n}  = \int_0^{\pi/2} \cos^{2n} x\, dx 
  = {(2n-1)!! \over (2n)!!}{\pi \over 2}
  = {\pi \over 2} \cdot \hbox{the coefficient of } x^{n} \hbox{ in } (1-x)^{-1/2} 
$$
i.e.,
$$
{1 \over \sqrt{1-x}}  = {2 \over \pi} \left(W_0+W_2x+W_4x^2+W_6x^3+\cdots\right)  
$$

\bigskip

\noindent
({\it iv}) 
Now
assume that $k=2n+1 \geq 1$ is odd.
In that case, we realize that    
\begin{align}
W_{2n+1} & =\int_0^{\pi/2} \cos^{2n+1} x\, dx = {(2n)!! \over (2n+1)!!}  \nonumber \\
   & =  \hbox{the reciprocal of the coefficient of } x^{n} \hbox{ in } (1-x)^{-3/2} \nonumber 
\end{align}
as a straightforward calculation yields 
$$
  (1-x)^{-3/2} = \sum_{n=0}^{\infty} {-3/2 \choose n} (-x)^n  
   = \sum_{n=0}^{\infty}  {(2n+1)!! \over (2n)!!}  x^n
$$
Thus,
$$
  (1-x)^{-3/2}  =  
    {1 \over W_1} + {1 \over W_3} x + {1 \over W_5} x^2 + {1 \over W_7} x^3 + \cdots
$$

\bigskip

\section{Summary}

In our exploration of a particular sawtooth wave we came across the family of definite integrals  
$I_k =\int_0^{\pi/2} x \cos^k x\,dx,$ for $k \geq 0.$
It turns out that these integrals, related to Wallis integrals $W_k=\int_0^{\pi/2} \cos^k x\,dx,$
have intriguing properties, some of which we present in this paper. Using an iterative formula
we compute $I_k,$ and then investigate a relationship between $I_k$ and the family
$J_k =\int_0^{\pi/2} x \cos^k 2x\,dx,$ for $k \geq 0.$
Simplifying one identity, we came across the expression ${2n \choose n}$
which is related to Catalan numbers.
We discover two interesting power series expansions involving $J_k$ and Wallis integrals.
In our view, this math-rich material could be incorporated into an advanced Calculus course,
or, broken-down into steps, given as a problem-solving activity.

\bigskip
 
\section{References}

\medskip

\newcommand{\referenceentry}[1]{\noindent\hangindent=15pt\hangafter=1\baselineskip 14pt #1}

\referenceentry{[1] Glaister, P. (2003). Factorial sums. {\it International Journal of Mathematical
Education in Science and Technology}, 34(2), 250-257.}\\ 
https://doi.org/10.1080/0020739031000158272

\medskip

\referenceentry{[2] Gould, H. $\&$ Quaintance, J. (2012). Double Fun with Double Factorials,
{\it Mathematics Magazine}, 85(3), 177-192.}\\ https://doi.org/10.4169/math.mag.85.3.177

\medskip

\referenceentry{[3] Stewart, J., Clegg D., $\&$ Watson, S. (2020). 
{\it Calculus: Early Transcendendals, Ninth Edition}. Brooks Cole.}
\end{document}